\documentclass[12pt]{amsart}
\usepackage{amssymb,amscd}
\usepackage[russian]{babel}
\usepackage{graphicx}
\usepackage{pdfpages}
\usepackage{placeins}
\newtheorem{theorem}{Теорема}[section]

\theoremstyle{definition}
\newtheorem{definition}{Определение}[section]

\theoremstyle{remark}

\ifx\pdfpagewidth\undefined \else
\fi
 \evensidemargin \oddsidemargin

\advance\hoffset by -1in
\advance\voffset by -1in

\usepackage[cp1251]{inputenc}

\author{И.~А.~Иванов-Погодаев, А.~Я.~Канель-Белов}

\thanks{Moscow Institute of Physics and Technology, Bar-Ilan University 
and College of Mathematics and Statistics, Shenzhen University, Shenzhen, 518061, China}

\address{Moscow Institute of Physics and Technology, Russia, Bar-Ilan University, Israel  }
\email{ivanov.pogodaev@mail.ru}
\address{College of Mathematics and Statistics, Shenzhen University, Shenzhen, 518061, China}
\email{kanel@mccme.ru}

\title{Конструкция бесконечной конечно определенной нильполугруппы}

\begin{document}

\let \mathbf=\texttt
\tabcolsep 2pt

\begin{abstract}

Работа посвящена конструкции конечно определенной бесконечной нильполугруппы, удовлетворяющей тождеству $x^9=0$. Эта конструкция отвечает на проблему  Л.~Н.~Шеврина и М.~В.~Сапира. 

Доказательство основано на построении последовательности геометрических комплексов, каждый из которых склеен из нескольких простых 4-циклов (квадратов). Комплексы обладают набором геометрических и комбинаторных свойств. Полугруппа строится как множество слов-кодировок путей на таких комплексах. Определяющие соотношения соответствуют парам эквивалентных коротких путей. Кратчайшим путям в смысле естественной метрики будут соответствовать ненулевые слова в полугруппе. Слова, не соответвующие путям на каком-либо комплексе или соответствующие некратчайшим путям, приводятся к нулю.
Полная версия работы --  \cite{nilsemigroup}.

Данная работа была проведена с помощью Российского Научного Фонда Грант N 17-11-01377.  Первый автор является победителем конкурса ``Молодая математика России''.

УДК 512.53, MSC: 20M05  

Ключевые слова: конечно определенные полугруппы, проблемы бернсайдовского типа; finitely presented semigroups, burnside-type problems.

\end{abstract}

\maketitle

\section{Введение} \label{nachalo}

Работа посвящена построению конечно определенных нильполугрупп. Доказана

\medskip
{\bf Теорема.} {\it Существует конечно определенная бесконечная нильполугруппа, удовлетворяющая тождеству $x^9=0$.}
\medskip

 Подобные построения затрагивают, с одной стороны, проблемы бернсайдовского типа, а с другой -- проблемы построения конечно определенных объектов. В полной версии работы  \cite{nilsemigroup} приведена более развернутая история вопроса.


Проблемы бернсайдовского типа внесли огромный вклад в развитие современной алгебры. Эта проблематика охватила большой круг вопросов,
как в теории групп, так и в смежных областях, стимулировала алгебраические исследования.

 Вопрос о локальной  конечности групп с тождеством $x^n = 1$ был решен отрицательно в знаменитых работах П.~С.~Новикова и С.~И.~Адяна  \cite{Novikov-Adyan}: После этого оценка улучшалась. Недавно С.~И.~Адян улучшил оценку до $n \ge 101$ (см.  \cite{Adyan1}).

 Работы П.~С.~Новикова и С.~И.~Адяна оказали огромное влияние на творчество И.~А.~Рипса, который в дальнейшем разработал метод канонической формы и построил примеры бесконечных периодических групп, обладающих дополнительными свойствами.

Все имеющиеся примеры бесконечных периодических групп бесконечно определены. Чрезвычайно глубоким и вдохновляющим является следующий открытый вопрос (входящий в список основных алгебраических проблем в теории групп):

\medskip

{\bf Вопрос.}\
{\it Существует ли конечно определенная бесконечная периодическая группа?
}
\medskip

Ответа на этот вопрос не известно как для случая, когда периоды ограничены в совокупности, так и для неограниченного случая.

\smallskip

На проблематику, связанную с построением разного рода
 экзотических объектов с помощью конечного числа определяющих соотношений обратил
внимание В.~Н.~Латышев. Им же была поставлена проблема существования конечно
определенного нилькольца.

\smallskip
{\bf Вопрос (В.~Н.~Латышев).} {\it Существует ли конечно определенное бесконечномерное нилькольцо?
}

\smallskip

Фундаментальную проблему существования конечно определенной нильполугруппы поставили Л.~Н.~Шеврин и М.~В.~Сапир в Свердловской Тетради (3.61б) \cite{Sverdlovsk},  вопрос 3.8 в \cite{Obzor}.

\medskip
{\bf Вопрос (Л.~Н.~Шеврин, М.~В.~Сапир).} {\it Существует ли конечно определенная бесконечная нильполугруппа?
}

\smallskip

{\bf Благодарности.} Авторы признательны руководителям семинара <<Теория колец>> на кафедре Высшей Алгебры механико-математического факультета МГУ В.~Н.~Латышеву и А.~В.~Михалеву за полезные обсуждения и постоянное, в течение ряда лет, внимание к работе. Мы также благодарны  И.~А.~Рипсу, Л.~Н.~Шеврину, А.~Х.~Шеню, Н.~К.~Верещагину, А.~Эршлер за полезные обсуждения, Ф.~Дюранду, Ц.~Селле, Л.~А.~Бокутю, Ю. ~Чэну, Т.~Фернику за поддержку в участии на конференциях,  А.~С.~Малистову за помощь в оформлении статьи. Особую благодарность мы выражаем А.~Л. ~Семенову за полезные советы и внимание к работе.

\section{План доказательства}      \label{ScFnDefInfGen}

\subsection{Схема построения.}

Пусть $W$~-- бесквадратное слово над алфавитом из трех букв. Если каждое его неподслово (т.е. антислово) объявить нулем, то естественно возникающая полугруппа слов обладает тождеством $x^2=0$. То есть построение бесконечной нильполугруппы с помощью бесконечной последовательности определяющий соотношений не представляет сложностей. Однако, в конечно определенном случае естественная конструкция, связанная с заданием множества нулевых слов как подслов слов из некоторого семейства работает плохо. Требуется привлечь новые идеи.

\smallskip

Существуют конечные наборы многоугольников (плиток) с заданными краевыми условиями, которыми можно замостить плоскость лишь непериодическим способом.
Широко известна, например, мозаика Пенроуза.
Можно использовать мозаику и пути на ней для введения определяющих соотношений в полугруппе. Буквы будут кодировать узлы стыков плиток и ребра между ними. В непериодической мозаике не будет периодических путей. Конечность числа краевых условий позволяет ввести конечное число определяющих соотношений.

\smallskip

Итак, основная идея -- построить граф, по которому будут проходить пути и использовать его как базу для введения определяющих соотношений.  Кодировки путей на этом графе будут элементами полугруппы.  Более точно, мы построим последовательность геометрических комплексов, составленных из простых 4-циклов ({\it плиток}).  Плиток в каждом комплексе конечное число, и прикладываются они одна к другой по целой стороне. Следующий комплекс последовательности получается из предыдущего с помощью иерархического правила разбиения плиток.  Полученная последовательность комплексов используется для введения определяющих соотношений в полугруппе.

\smallskip
Имея построенную последовательность комплексов, далее будем действовать по следующему плану:

1. Все комплексы последовательности в совокупности содержат конечное число конфигураций вершин -- узлов, где сходятся несколько плиток.  Выпишем все возможные типы таких узлов и обозначим их буквами алфавита. 

2. Вторую серию букв алфавита введем для всевозможных входящих и выходящих ребер во все вершины. 

3. Теперь последовательность букв (слово) может кодировать последовательность вершин, которые мы проходим вдоль пути, лежащем на произвольном комплексе. Выпишем все последовательности букв длины $7$ которые не представляют никакого пути ни на одном из комплексов.  Занесем все такие последовательности в список запрещенных и объявим их нулями введя соответствующие определяющие соотношения. 

4. Рассмотрим произвольную плитку $T$. Существует два кратчайших пути, соединяющих ее противоположные вершины -- по двум соседним сторонам или по двум другим соседним сторонам. Такую пару путей объявим {\it эквивалентными}. Проделаем эту операцию для всех возможных комбинаций 4 цветов вершин, составляющих плитку. Для каждой плитки вводятся две пары эквивалентных путей. 
Всего мы вводим конечное число таких эквивалентностей, так как существует конечное множество типов плиток.  Отметим, что последовательности, кодирующая указанные пути, имеет длину $7$: вершина $A$, выходящее ребро, входящее ребро, вершина $B$, выходящее ребро, входящее ребро, вершина $C$.  Для каждой такой пары эквивалентных путей мы вводим соответствующее определяющее соотношение для кодирующих их слов.

5. Также внесем в список запрещенных пути вида туда-и-обратно по одному ребру. Для каждого пути вида туда-и-обратно выпишем его кодировку и соответствующую последовательность букв объявим нулем, введя определяющее соотношение. Все такие соотношения будут иметь длину $7$, их также конечное число.

6. Таким образом, с помощью локальных замен эквивалентных кусков на эквивалентные, можно переводить одни пути в другие. {\it Запрещенными} мы будем считать пути, кодировки которых непосредственно объявлены нулями в определяющих соотношениях. {\it Нулевыми} путями будем считать пути, кодировки которых приводятся локальными заменами к форме, содержащей запрещенный подпуть.

7. Заметим, что все введенные соотношения имеют длину $7$. В частности, их конечное число. Кроме того, эквивалентные переходы не меняют длины пути. Далее мы проверяем, что кратчайший путь в любом комплексе не может быть приведен к виду содержащему запрещенный участок. Также мы проверяем, что некратчайший путь может быть приведен локальными заменами к виду, содержащему путь туда-и-обратно по одному ребру.

8. Далее, мы проверяем,  что любой путь, содержащий периодический участок, может быть приведен к виду, содержащему запрещенный участок. То есть, его кодировка приводится к нулю с помощью определяющих соотношений.

9. Комплексы последовательности обладают самоподобием, комплекс более высокого уровня содержит и пути комплексов более низких уровней. То есть, локально комплексы более высоких уровней устроены также, как комплексы более низких уровней. Фактически, они состоят из одних и тех же плиток. Для всех комплексов  используется одно конечное множество определяющих соотношений, ограниченных по длине.

10. Таким образом, с последовательностью комплексов будет связана конечно определенная полугруппа, где словам соответствуют пути. Слова $A$ и $B$ равны в полугруппе, если существует цепочка локальных эквивалентостей слов $A=A_1\equiv A_2\equiv \dots A_k=B$, такая что для каждого $i<k$ можно выбрать комплекс из последовательности и два пути $C$ и $D$ на нем, с одинаковыми началом и концом, причем: 

1) кодировка $C$ соответствует слову $A_i$, кодировка $D$ соответствует $A_{i+1}$;

2) Путь $D$ получается из $C$ заменой одного эквивалентного участка на другой.

 Такая полугруппа будет конечно определенной и бесконечной нильполугруппой.

\medskip

\subsection{Свойства базового геометрического комплекса}

Теперь обсудим свойства комплексов, входящих в последовательность.


\medskip
\begin{description}
  \item[1. Локальная конечность] Каждый комплекс семейства состоит из конечного числа 4-циклов (плиток), соединенных друг с другом по сторонам. Узлы -- это вершины плиток, существует конечное число типов узлов, и конечное число типов выходящих из них ребер.
  
  \item [2. Равномерная эллиптичность] Любой достаточно длинный путь, соединяющий узлы $A$ и $B$ может быть переведен локальными заменами в другой путь, отличающийся достаточно сильно от начального.

      Более точно, пусть расстояние между двумя путями длины $n$ с общими концами -- это максимум по расстояниям между их соответствющими точками.  Максимум такого расстояния для выбранных двух концов пути будем называть {\ шириной пучка}.  Равномерная эллиптичность комплекса означает что для любого $k$ можно выбрать $n$ так, что для любых точек на расстоянии $n$ ширина соответстующего пучка будет более $k$.
      
   \item [3. Детерминированность] Можно выбрать глобальную константу $k$ и покрасить узлы всех комплексов последовательности в $k$ цветов так, чтобы цвета трех вершин любой плитки однозначно определяли цвет четвертой вершины этой плитки.
   
  \item [4. Апериодичность] На комплексах не должно быть путей, отвечающих периодическим словам.
\end{description}

\smallskip

Комплекс с такими свойствами не может содержать в себе кратчайших путей, выражаемых периодическим словом. После введения определяющих соотношений по плану, указанному выше, мы получаем возможность локальных преобразований слов в полугруппе.

 Иерархическая структура и равномерная эллиптичность комплексов позволяют применить редукционный процесс к слову. 
 Для любого слова есть две возможности: либо оно приводится к канонической форме и мы получаем, что оно является кодировкой пути на некотором комплексе, либо в процессе редукции возникает запрещенное слово, и мы приходим к нулю.
 Кодировки кратчайших пути, лежащих на каком-либо комплексе, не приводятся к нулю. Для любой заданной длины $n$ в последовательности комплексов можно выбрать достаточно большой комплекс, содержащий кратчайшие пути длины $n$.
Таким образом, полугруппа, соответствующая множеству кодировок кратчайших путей на построенной последовательности комплексов, будет бесконечной конечно определенной нильполугруппой.

\subsection{Основные составные части доказательства}

\begin{figure}[hbtp]
\centering
\includegraphics[width=0.6\textwidth]{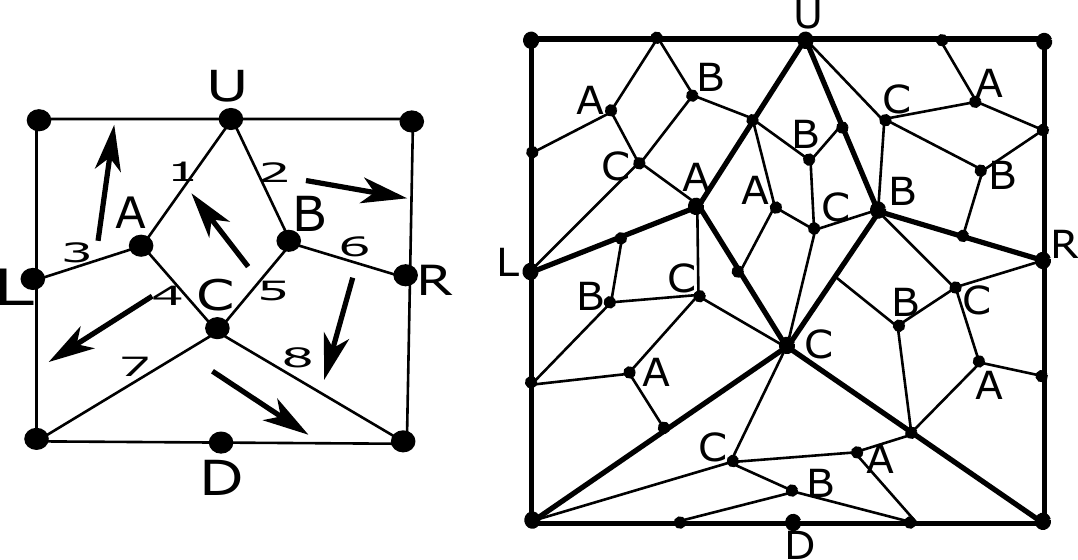}
\caption{Схема разбиения и макроплитка третьего уровня с типами вершин }
\label{fig:macro}
\end{figure}

\begin{definition}
{\it Макроплиткой}  $T_n$ уровня $n$ будем называть:
 
 Для $n=1$ -- простой цикл, состоящий из $4$ вершин и $4$ ребер, среди ребер мы выделяем верхнее, правое, левое и нижнее, в соответствии с естественной ориентацией.
 
 Для $n>1$ -- макроплитку $T_{n-1}$, к которой применена операция {\it разбиения}: каждая ориентированная макроплитка $1$ уровня, из которой состоит $T_{n-1}$, разбивается на шесть макроплиток $1$ уровня, ориентированных в соответствии со схемой на левой стороне рисунка~\ref{fig:macro}.  Стрелками показаны верхние ребра в каждой из шести макроплиток разбиения.  При разбиении образуется $7$ вершин. Три из них -- внутри первоначальной плитки, будем говорить, что они имеют типы $A$, $B$ и $C$. Еще четыре образуются на серединах верхней, правой, нижней и левой сторон: $U$, $R$, $D$, $L$ соответственно.
\end{definition}

 В соответствии со схемой иерархического разбиения, при переходе от $T_n$ к $T_{n+1}$ добавляется несколько вершин. Некоторые вершины являются общими для соседних макроплиток, прилегающих друг к другу по стороне. Все вершины в макроплитке, кроме изначальных четырех (углов) создаются при каком-либо из разбиений.

Таким образом, макроплитка $T_n$ состоит из $6^{n-1}$ макроплиток $1$ уровня, прилегающих друг к другу в соответствии с иерархической схемой.

\begin{definition}
Проводимые при разбиении ребра будем называть {\it внутренними ребрами}, принадлежащими разбиваемой макроплитке. {\it Уровнем ребра} будем называть уровень макроплитки, которой оно принадлежит. 
Будем считать, что  изначальные четыре вершины (углы) имеют {\it глубину} равную $-1$. Создаваемые при первом разбиении вершины $A$, $B$, $C$, $U$, $L$, $R$, $D$ имеют {\it глубину} 0. Определим глубину новых вершин. Все создаваемые при разбиении вершины получают глубину на 1 больше, чем максимальная глубина вершины до разбиения.
Из построения ясно, что у любого ребра макроплитки хотя бы один из концов является вершиной, созданной на предыдущем шаге. Следовательно, создаваемая вершина в середине ребра получает глубину на $1$ большую, чем, максимальная глубина двух концов ребра.
{\it Типами} внутренних ребер будем называть $8$ видов ребер образующихся при разбиениях, они отмечены как $1-8$ на левой части рисунка~\ref{fig:macro}. Есть также $4$ типа граничных ребер, левое, правое, верхнее и нижнее. Граничные ребра комплекса относятся к этим типам. К ним же относятся границы  подклееных макроплиток (см ниже).

\end{definition}

\begin{definition}
Определим операцию {\it подклейки}. 

Мы рассматриваем все пути $X_1X_2YZ_2Z_1$ из пяти  вершин (и четырех ребер), такие что:

1) Вершины $X_1$, $Y$, $Z_1$ {\bf не} являются тремя углами из четырех никакой макроплитки;

2) Вершины $X_1$, $Z_1$ являются боковыми вершинами глубины $k-1$, где $k$ -- максимальная глубина вершины;

3) Вершины $X_2$ является серединой $X_1Y$, а $Z_2$ -- серединой $Z_1Y$, то есть,  $X_2$ и $Z_2$  являются боковыми вершинами глубины $k$, созданными  самыми последними, при последней операции  разбиения.

4) Вершина $Y$ имеет глубину $k-2$.

5) Путь $X_1X_2YZ_2Z_1$ выбран так, что уровень ребра, на котором лежит вершина $X_1$ больше уровня ребра вершины $Z_1$.

В случае, когда обе вершины $X_1$ и $Z_1$ имеют один уровень ребра, ориентация пути выбирается следующим образом.

случай А) $X_1$ и $Z_1$ лежат на ребрах разного типа. Тогда $X_1$ выбирается как вершина
с большим номером типа ребра.

случай B) $X_1$ и $Z_1$ лежат на одном ребре. Тогда для внутреннего ребра $X_1$ выбирается как вершина лежащая ближе к краю макроплитки. Для внешнего ребра $X_1$ выбирается как предшествующая $Z_1$ при обходе контура по часовой стрелке. 

\medskip

Далее, для каждого такого пути создаются шесть новых вершин $T_1$, $T_2$, $T_3$, $T_A$,  $T_B$,  $T_C$, не лежащих в плоскости $X_1YZ_1$, а также проводятся новые ребра $X_1 T_2$, $X_2 T_A$, $X_2 T_B$, $T_2 T_B$, $T_C T_B$, $T_C T_A$, $T_2 T_1$, $T_3 T_1$, $T_3 Z_1$, $T_C Z_1$, $T_A Z_2$  и появляется новая {\it подклееная } макроплитка $X_1YZ_1T_1$ (рисунок~\ref{fig:pasting}).

\begin{figure}[hbtp]
\centering
\includegraphics[width=0.4\textwidth]{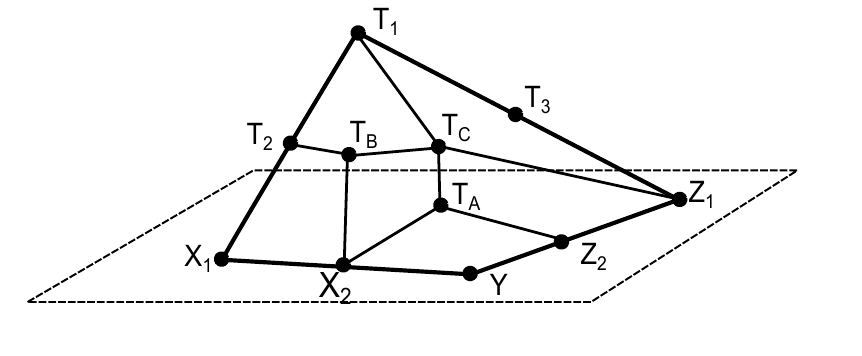}
\caption{Подклейка.}
\label{fig:pasting}
\end{figure}

Созданные вершины будем называть {\it подклеенными}, причем $T_1$ присваивается глубина $k-1$, а вершинам $T_2$, $T_3$, $T_A$,  $T_B$,  $T_C$ -- глубина $k$.
Вершину $T_1$ будем считать угловой, вершины $T_2$, $T_3$ -- боковыми, а вершины $T_A$,  $T_B$,  $T_C$ -- внутренними.

Тип вершин вдоль пути $X_1X_2YZ_2Z_1$ не меняется при проведении подклейки.

Макроплитку, которой принадлежит вершина $Y$, будем называть {\it базовой плоскостью} подклейки.

\medskip

{\bf Примечание 1.} Заметим, что в подклееной макроплитке по построению определяется верхняя, а, следовательно, и остальные стороны: она выглядит так же, как если бы к плитке $X_1YZ_1T_1$ применили разбиение, считая сторону $X_1Y$ верхней. Также можно считать, что в этой макроплитке $T_2$ середина стороны $T_1X_1$, а $T_3$ -- середина стороны $T_1Z_1$.

\medskip

\end{definition}

\begin{definition}
{\it Комплексом}  $K_n$ уровня $n$ будем называть:
 
 Для $n=1,2,3$ -- макроплитку $T_n$;
 
 Для $n>3$ -- комплекс $T_{n-1}$, к которому сначала применена операция разбиения, а затем -- операция подклейки. 
\end{definition}

{\bf Замечание 1.} Начиная с четвертого уровня комплекса (и при построении пятого уровня) путь $X_1YZ_1$ из определения операции подклейки может частично лежать в базовой плоскости, а частично в подклееной части (если $Z_1$ появилась при предыдущей подклейке в качестве вершины $T_2$ или $T_3$). На больших уровнях комплекса весь путь из определения подклейки может лежать полностью в подклееных плитках. Таким образом, возникают подклейки к частям подклееных когда-то макроплиток и т.д.

{\bf Замечание 2.} В итоге мы всегда будем рассматривать комплексы конечного уровня. Предельного перехода применяться не будет.

\begin{definition}
Пусть $d(A,B)$ -- стандартное расстояние  между точками $A$ и $B$, длина пути по ребрам. {\it Пучком путей $\Omega_n(A,B)$ с концами в точках $A$ и $B$} будем называть множество кратчайших путей на комплексе $K_n$, соединяющих точки $A$ и $B$.  Пусть $s$ --  длина всех этих путей. Тогда для  пути $P$ определены точки комплекса, по которым он проходит: $P_0=A$, $P_1$, \dots , $P_{s-1}$, $P_s=B$. На $\Omega_n(A,B)$ можно определить расстояние между путями. Пусть $P$, $Q\in \Omega_n(A,B)$. Расстоянием между путями $r(P,Q)$ будем называть $\max_{0\le i\le s}{d(P_i,Q_i)}$, то есть максимальное расстояние между соответствующими точками двух путей с одинаковыми началом и концом. 

{\it Шириной} пучка $\Omega_n(A,B)$ будем называть максимум $r(P,Q)$, где  $P$, $Q\in \Omega_n(A,B)$, то есть  максимальное расстояние между принадлежащими ему двумя путями.

\end{definition}

\begin{theorem}[О существовании] \label{complex}

В построенной последовательности геометрических комплексов $K_n$ выполнены следующие свойства:

1. Локальная конечность. Существует конечное число типов вершин, а также типов входящих и выходящих ребер.

2. Равномерная эллиптичность. Для любого натурального числа $k$ существуют натуральные числа $N$ и $S$, такие что выполнено следующее свойство: в каждом комплексе $K_n$ для $n>N$ для любых точек $A$ и $B$ c расстоянием не менее $S$, пучок путей $\Omega_n(A,B)$ имеет ширину не менее $k$.
\end{theorem}

\begin{theorem}[О детерминированности] \label{determin}

Построенная последовательность комплексов $K_n$ обладает свойством детерминированности:

cуществует такая глобальная константа $N$, что все точки каждого из комплексов $K_n$ можно покрасить в $N$ цветов так, чтобы для любых четырех точек, образующих макроплитку $1$ уровня в комплексе $K_n$ по цветам трех из них можно однозначно установить цвет четвертой.
\end{theorem}

C построенной последовательностью комплексов $K_n$ можно связать конечно определенную полугруппу $S$ следующим образом.

Буквам полугруппы отвечают различные цвета для вершин комплекса, а также различные типы входящих и выходящих в вершины ребер. Словам отвечают последовательности букв, то есть кодировки путей на комплексе. 

Для слов $w$ длины $<10$, не кодирующих никакой путь на комплексе, а также для слов, кодирующих путь туда-обратно по одному и тому же ребру, вводятся определяющие соотношения $w=0$. 

Для каждой макроплитки уровня $1$ семейства комплексов вводится определяющее соотношение $w=v$, где $w$ --  кодировка пути по двум соседним сторонам макроплитки, а $v$ -- кодировка пути с теми же началом и концом, но по двум другим сторонам. Кодировка каждого такого пути включает $7$ букв -- цвет начальной вершины, тип ребра выхода, тип ребра входа во вторую вершину, цвет второй вершины, тип ребра выхода из второй вершины, тип ребра входа в конечную вершину, цвет конечной вершины.

\begin{theorem}[О полугруппе, связанной с последовательностью комплексов] \label{semigroup}

Полугруппа $S$, связанная с последовательностью $K_n$ обладает следующими свойствами:

Слова, отвечающие путям различной длины, различны в полугруппе;

Слово, отвечающее кодировке пути по стороне макроплитки любого уровня не может быть приведено к нулю определяющими соотношениями;

Слово, содержащее девятую степень некоторого подслова может быть приведено к нулю с помощью определяющих соотношений.

\end{theorem}

\end{document}